\documentclass[10pt]{article}
\usepackage{latexsym}
\usepackage{amsbsy}
\usepackage{amsfonts,mathrsfs,amssymb}
\usepackage{amsmath,amsthm}
\usepackage[latin1]{inputenc}
\usepackage{graphics,graphicx}
\usepackage{a4wide}

 \newtheorem{thmm}{Theorem}
 \newtheorem{defin}[thmm]{Definition}
 \newtheorem{lemm}[thmm]{Lemma}
 \newtheorem{propp}[thmm]{Proposition}
 \newtheorem{corr}[thmm]{Corollary}
 \newtheorem{remm}[thmm]{Remark}
 \newtheorem{exx}[thmm]{Example}

 \newcommand{\bthm}{\begin{thmm}}
 \newcommand{\ethm}{\end{thmm}}
 \newcommand{\bd}{\begin{defin}}
 \newcommand{\ed}{\end{defin}}
 \newcommand{\blem}{\begin{lemm}}
 \newcommand{\elem}{\end{lemm}}
 \newcommand{\bcor}{\begin{corr}}
 \newcommand{\ecor}{\end{corr}}
 \newcommand{\bprop}{\begin{propp}}
 \newcommand{\eprop}{\end{propp}}
 \newcommand{\brem}{\begin{remm} \rm}
 \newcommand{\erem}{\end{remm}}
 \newcommand{\bex}{\begin{exx} \rm}
 \newcommand{\eex}{\end{exx}}

 \newcommand{\pr}{\noindent{\bf Proof. }}
 \newcommand{\ep}{\nolinebreak{\hspace*{\fill}$\Box$ \vspace*{0.25cm}}}

 \newcommand{\beq}{\begin{equation}}
 \newcommand{\eeq}{\end{equation} }

 \newcommand{\bea}{\begin{eqnarray}}
 \newcommand{\eea}{\end{eqnarray}}

 \newcommand{\beas}{\begin{eqnarray*}}
 \newcommand{\eeas}{\end{eqnarray*}}

 \newcommand{\beqs}{\begin{equation*}}
 \newcommand{\eeqs}{\end{equation*}}

 \newcommand{\bi}{\begin{itemize}}
 \newcommand{\ei}{\end{itemize}}

 \newcommand{\ben}{\begin{enumerate}}
 \newcommand{\een}{\end{enumerate}}

 \newcommand{\ba}{\begin{array}}
 \newcommand{\ea}{\end{array}}

 \newcommand{\R}{\mathbb R}
 \newcommand{\N}{\mathbb N}

 \newcommand{\cA}{\ensuremath{{\cal A}}}
 \newcommand{\cC}{\ensuremath{{\cal C}}}

 \newcommand{\cL}{\ensuremath{{\cal L}}}

 \newcommand{\pd}{\partial}

 \newcommand{\vphi}{\varphi}

 \newcommand{\vf}{{\bf v}}

\begin{document}

 \title{
 Variational problems with fractional derivatives:
 Invariance conditions and N\"{o}ther's theorem
 }

 \author{
 Teodor M. Atanackovi\'c
 \footnote{Faculty of Technical Sciences, Institute of Mechanics, University of Novi Sad, Trg Dositeja Obradovi\' ca 6, 21000 Novi Sad, Serbia.
 Electronic mail: atanackovic@uns.ac.rs.
 Work supported by
 Projects 144016 and 144019 of the Serbian Ministry of Science
 and
 START-project Y-237 of the Austrian Science Fund.}\\
 Sanja Konjik
 \footnote{Faculty of Agriculture, Department of Agricultural Engineering, University of Novi Sad, Trg Dositeja Obradovi\' ca 8, 21000 Novi Sad, Serbia.
 Electronic mail: sanja\_konjik@uns.ac.rs. }\\
 Stevan Pilipovi\'c
 \footnote{Faculty of Sciences, Department of Mathematics, University of Novi Sad, Trg Dositeja Obradovi\' ca 4, 21000 Novi Sad, Serbia.
 Electronic mail: pilipovic@dmi.uns.ac.rs}\\
 Srboljub Simi\'c
 \footnote{Faculty of Technical Sciences, Institute of Mechanics, University of Novi Sad, Trg Dositeja Obradovi\' ca 6, 21000 Novi Sad, Serbia.
 Electronic mail: ssimic@uns.ac.rs}\\
 }

 \date{}
 \maketitle

 \begin{abstract}
 A variational principle for Lagrangian densities containing
 derivatives of real order is formulated and the invariance of
 this principle is studied in two characteristic cases. Necessary
 and sufficient conditions for an infinitesimal transformation
 group (basic N\"{o}ther's identity) are obtained. These conditions
 extend the classical results, valid for integer order derivatives.
 A generalization of N\"{o}ther's theorem leading to conservation
 laws for fractional Euler-Lagrangian equation is obtained as well.
 Results are illustrated by several concrete examples. Finally, an
 approximation of a fractional Euler-Lagrangian equation by a system
 of integer order equations is used for the formulation of an
 approximated invariance condition and corresponding conservation laws.

 \vskip5pt

 \noindent {\bf Mathematics Subject Classification (2000):}
 Primary: 49K05; secondary: 26A33

 \vskip5pt

 \noindent {\bf PACS numbers:}
 02.30 Xx, 45.10 Hj

 \vskip5pt

 \noindent {\bf Keywords:}
 variational problem, Riemann-Liouville fractional derivatives,
 variational symmetry, infinitesimal criterion, N\"{o}ther's theorem,
 conservation laws, approximations
 \end{abstract}

 \section{Introduction}

 There are two distinct approaches in the formulation of fractional
 differential equations of models in various branches of science, e.g.
 physics. In the first one, differential
 equations containing integer order derivatives are modified by
 replacing one or more of them with fractional ones
 (derivatives of real order). In the second approach, which we will
 follow in the sequel, one starts with a variational formulation
 of a physical process in which one modifies the Lagrangian density by
 replacing integer order derivatives with fractional ones. Then
 the action integral in the sense of Hamilton is minimized and
 the governing equation of a physical process is obtained. Hence, one
 is faced with the following problem: find minima (or maxima) of a
 functional
 \beq \label{eq:prva}
 \cL[u]=\int_A^B L(t,u(t),{}_aD_t^\alpha u)\, dt,
 \quad 0< \alpha< 1,
 \eeq
 where ${}_aD_t^\alpha u$ is the left Riemann-Liouville fractional
 derivative, under certain assumptions on a Lagrangian
 $L$, as well as on functions $u$ among which minimizers are sought.
 After this, the Euler-Lagrange equations are formed for the modified
 Lagrangian, which leads to equations that intrinsically characterize
 a physical process. This approach has a more sound physical basis (see
 e.g. \cite{Bal06, Bal04, Bal05, Rekh, SabatierAgrawalTeMa}).
 It provides the possibility of formulating conservation laws via
 N\"{o}ther's theorem. Conservation laws are very important in
 particle reaction physics for example, where they are often
 postulated as conservation principles \cite{Schulte}.
 The central point is the fact that equations of a model, i.e.,
 Euler-Lagrange equations, give minimizers of a functional.
 As a rule, fractional differential equations obtained through
 Euler-Lagrange equations contain left {\em and} right fractional
 derivatives, which makes this approach more delicate.

 In the last years fractional calculus has become popular as a
 useful tool for solving problems from various fields
 \cite{Hilfer, WestBolognaGrigolini}
 (see also \cite{APZ07, GorenfloAbdel, GorenfloLuchkoMainardi,
 GorenfloMainardiVivoli-1, GorenfloMainardiVivoli-2,
 Jumarie01, Jumarie06, Lazopoulos06, MainardiPagniniGorenflo,
 Sakakibara, Stanislavsky}).

 The study of fractional variational problems also has a long
 history. F. Riewe \cite{Riewe96, Riewe97}
 investigated
 nonconservative Lagrangian and Hamiltonian mechanics and
 for those cases formulated a version of the Euler-Lagrange
 equations. O. P. Agrawal continued
 the study of the fractional Euler-Lagrange equations
 \cite{Agrawal02, Agrawal06, Agrawal07}, for general fractional
 variational problems involving
 Riemann-Liouville, Caputo and Riesz fractional derivatives.
 G. S. F. Frederico
 and D. F. M. Torres
 \cite{FredericoTorres} studied invariance properties of fractional
 variational problems with N\"{o}ther type theorems by
 introducing a new concept of fractional conserved quantity
 which is not constant in time, so the term conserved quantity is
 not clear.
 Conservation laws and Hamiltonian type equations for the
 fractional action principles have also been derived in
 \cite{TarasovZaslavsky}.
 There are nowadays numerous applications of fractional
 variational calculus, see e.g.
 \cite{AOP07, AtanackovicStankovic07, Bal06, Bal04, Bal05,
 DreisigmeyerYoung04, DreisigmeyerYoung03, GaiesElAkrmi,
 Klimek01, Lazopoulos06, Odibat07, Odibat06, Baleanu&co07}.

 We can summarize the novelties of our paper as follows. First, we
 derived infinitesimal criterion for a local one-parameter group
 of transformations to be a variational symmetry group for the
 fractional variational problem (\ref{eq:prva}). In previous work
 \cite{FredericoTorres, FredericoTorres-cl, FredericoTorres-C}
 this was done only in a special case. Moreover, we separately consider
 cases when the lower bound $a$ in the left Riemann-Liouville
 fractional derivative is not transformed (Subsection
 \ref{subsec:a not transformed}), and when $a$ is transformed
 (Subsection \ref{subsec:a transformed}). Both cases have
 their physical interpretation.
 In the case when fractional derivatives model memory effects the
 lower bound $a$ in the definition of derivative should not be
 varied. However, if fractional derivatives model nonlocal
 interactions (e.g. in nonlocal elasticity) the lower bound has to
 be varied.
 Second, we give a N\"{o}ther type theorem in terms of
 conservation laws as it is done in the classical theory. This
 approach preserves the essential property of conserved
 quantities to be constant in time.
 Third novelty is the approximation procedure which is given for
 the Euler-Lagrange equations and infinitesimal criterion,
 therefore also for N\"{o}ther's theorem. By the use of additional
 assumptions we obtain that appropriate sequences of classical
 Euler-Lagrange equations, infinitesimal criteria and conservation
 laws converge in the sense of an appropriate space of generalized
 functions to the corresponding Euler-Lagrange equations,
 infinitesimal criterion and conservation laws of the fractional
 variational problem (\ref{eq:prva}).
 There are also a number of illustrative examples which complete
 the theory.

 The paper is organized as follows.
 To the end of this introductory
 section we provide basic notions and definitions from the
 calculus with derivatives and integrals of any real order, and
 recall the Euler-Lagrange equations of (\ref{eq:prva}) for
 $(A,B)\subseteq (a,b)$ (so far only the case $(A,B)=(a,b)$ was
 treated, cf.\ \cite{Agrawal02, Agrawal06}).
 Section \ref{sec:inf inv}
 is devoted to the study of variational symmetry groups and invariance
 conditions.
 In Section \ref{sec:Noether} we derive a version of N\"{o}ther's
 theorem which establishes a relation between fractional
 variational symmetries and conservation laws for the Euler-Lagrange
 equation, which have the
 important property of being constant in time.
 We use in Section \ref{sec:approximations} the approximation of
 the Riemann-Liouville fractional derivative by a finite sum of
 classical derivatives. In this way and by introducing a suitable
 space of analytic functions we obtain approximated
 Euler-Lagrange equations, approximated local Lie group actions
 with corresponding infinitesimal criteria and N\"{o}ther
 theorems.

 \subsection{Preliminaries}

 In the sequel we briefly recall some basic facts from the
 fractional calculus.

 For $u\in L^1([a,b])$, $\alpha>0$ and $t\in [a,b]$ we define the left,
 resp.\ right Riemann-Liouville fractional integral of order $\alpha$,
 as
 $$
 {}_aI_t^\alpha u=\frac{1}{\Gamma (\alpha)}
 \int_a^t (t-\theta)^{\alpha -1} u(\theta)\, d\theta,
 \quad
 \mbox{ resp. }
 \quad
 {}_tI_b^\alpha u=\frac{1}{\Gamma (\alpha)}
 \int_t^b (\theta-t)^{\alpha -1} u(\theta)\, d\theta.
 $$
 Left, resp.\ right Riemann-Liouville fractional derivative
 of order $\alpha$, $0\leq\alpha<1$, is well defined for an
 absolutely continuous function $u$ in $[a,b]$, i.e. $u\in
 AC([a,b])$, and $t\in [a,b]$ as
 \beq \label{eq:RL-left}
 {}_aD_t^\alpha u=\frac{d}{dt} {}_aI_t^{1-\alpha} u=
 \frac{1}{\Gamma (1-\alpha)}\frac{d}{dt}
 \int_a^t \frac{u(\theta)}{(t-\theta)^{\alpha}}\, d\theta,
 \eeq
 resp.
 $$
 {}_tD_b^\alpha u=\Big(-\frac{d}{dt}\Big) {}_tI_b^{1-\alpha} u=
 \frac{1}{\Gamma (1-\alpha)}\Big(-\frac{d}{dt}\Big)
 \int_t^b \frac{u(\theta)}{(\theta-t)^{\alpha}}\, d\theta.
 $$
 If $f,g\in AC([a,b])$ and $0\leq\alpha<1$ the following fractional
 integration by parts formula holds:
 \beq \label{eq:frac int by parts}
 \int_a^b f(t) {}_aD_t^\alpha g\, dt = \int_a^b g(t)
 {}_tD_b^\alpha f\, dt.
 \eeq

 For the left Riemann-Liouville fractional derivative
 (\ref{eq:RL-left}) of order
 $\alpha$ ($0\leq\alpha<1$) we have:
 \beq \label{eq:relation left rl c}
 {}_aD_t^\alpha u= \frac{1}{\Gamma (1-\alpha)}
 \int_a^t \frac{\dot{u}(\theta)}{(t-\theta)^{\alpha}}\, d\theta
 + \frac{1}{\Gamma (1-\alpha)}
 \frac{u(a)}{(t-a)^\alpha},
 \quad t\in [a,b],
 \eeq
 The integral on the right hand side is called the left
 Caputo fractional derivative of order $\alpha$ and is denoted
 by ${}_a^cD_t^\alpha u$. Similarly, the right Caputo fractional
 derivative ${}_t^cD_b^\alpha u$ is defined.

 It follows from (\ref{eq:relation left rl c}) that
 the left Riemann-Liouville equals to the left Caputo fractional
 derivative in the case $u(a)=0$ (the analogue holds for the right
 derivatives under the assumption $u(b)=0$). The same condition,
 i.e. $u(a)=0$, provides that $\frac{d}{dt}{}_aD_t^\alpha u=
 {}_aD_t^\alpha \frac{d}{dt} u$.

 We introduce the notation which will be
 used throughout this paper: derivatives of
 Lagrangian $L=L(t,u(t), {}_aD_t^\alpha u)$ with
 respect to the first, second and third variable will be denoted
 by $\displaystyle \frac{\pd L}{\pd t}$, $\displaystyle \frac{\pd L}{\pd u}$ and
 $\displaystyle \frac{\pd L}{\pd {}_aD_t^\alpha u}$, or by $\pd_1 L$, $\pd_2 L$
 and $\pd_3 L$ respectively.

 \subsection{Euler-Lagrange equations}

 As stated in Introduction, we solve a fractional variational
 problem
 \beq \label{eq:fractional vp}
 \cL[u]=\int_A^B L(t,u(t),{}_aD_t^\alpha u)\, dt \to \min,
 \quad 0< \alpha< 1.
 \eeq
 $(A,B)$ is a subinterval of $(a,b)$, and Lagrangian $L$
 is a function in $(a,b)\times \R\times \R$ such that
 \beq \label{eq:conditions on L}
 \left.\ba{c}
 L\in\cC^1((a,b)\times \R\times \R)\\
 \mbox{ and}\\
 t\mapsto\pd_3 L(t,u(t),{}_aD_t^\alpha u)
 \in AC([a,b]), \mbox{ for every } u\in AC([a,b])
 \ea\right\}
 \eeq
 Solutions to (\ref{eq:fractional vp}) are sought among all
 absolutely continuous functions in $[a,b]$, which
 in addition satisfy condition $u(a)=a_0$, for a fixed
 $a_0\in\R$.

 One can consider more general problems with Lagrangians depending
 also on the right Riemann-Liouville fractional derivative. Our
 results can be easily formulated in that case, and because of
 simplicity we shall consider only Lagrangians of the form
 (\ref{eq:fractional vp}). Moreover, in (\ref{eq:fractional vp})
 one can consider the case $\alpha>1$, but this also does not give
 any essential novelty and because of that it is skipped.

 We will recall results about the Euler-Lagrange equations
 obtained in \cite{Agrawal02, Agrawal06, fractionalEL}.
 As we said, we consider the fractional
 variational problem defined by (\ref{eq:fractional vp}).

 Fractional Euler-Lagrange equations, which provide a
 necessary condition for extremals of a fractional variational
 problem, have been recently studied in
 \cite{Agrawal02, Agrawal06, fractionalEL}.
 Let $A=a$ and $B=b$. Euler-Lagrange equations
 are obtained in \cite{Agrawal02, Agrawal06}:
 \beq \label{eq:EL with RL}
 \frac{\pd L}{\pd u} + {}_tD_b^\alpha \Big(\frac{\pd L}{\pd _aD_t^\alpha
 u}\Big) =0.
 \eeq
 In terms of the Caputo fractional derivative,
 the Euler-Lagrange equation (\ref{eq:EL with RL})
 is given in \cite{fractionalEL} and reads:
 \beq \label{eq:EL popravljene}
 \frac{\pd L}{\pd u} + {}_t^cD_b^\alpha \Big(\frac{\pd L}{\pd _aD_t^\alpha
 u}\Big) +\frac{\pd L}{\pd {}_aD_t^\alpha u}{\Big\vert}_{_{t=b}}
 \frac{1}{\Gamma(1-\alpha)}\frac{1}{(b-t)^\alpha}=0,
 \eeq

 Euler-Lagrange equations for (\ref{eq:fractional vp})
 are derived in \cite{fractionalEL}:
 \bea
 \frac{\pd L}{\pd u}+{}_{t}^cD_{B}^{\alpha}
 \Big( \frac{\pd L}{\pd {}_{a}D_{t}^{\alpha}u}\Big)
 +\frac{\pd L}{\pd {}_aD_t^\alpha u}{\Big\vert}_{_{t=B}}
 \frac{1}{\Gamma(1-\alpha)}\frac{1}{(B-t)^\alpha}
 &=& 0,\; t\in (A,B) \label{eq:EL AB}\\
 {}_{t}D_{B}^{\alpha}\Big( \frac{\pd L}{\pd {}_{a}D_{t}^{\alpha}u}
 \Big) -{}_{t}D_{A}^{\alpha}\Big(
 \frac{\pd L}{\pd {}_{a}D_{t}^{\alpha}u}\Big)  &=& 0,\; t\in
 (a,A). \label{eq:EL aA}
 \eea
 Equation (\ref{eq:EL AB}) is equivalent to
 \beq \label{eq:EL AB RL}
 \frac{\pd L}{\pd u} + {}_tD_B^\alpha \Big(\frac{\pd L}{\pd _aD_t^\alpha
 u}\Big) =0, \quad t\in (A,B).
 \eeq

 \brem \label{rem:1der1}
 If we replace the Lagrangian in (\ref{eq:fractional vp}) by $L(t,
 u(t), {}_aD_t^\alpha u, \dot{u}(t))$ then Euler-Lagrange equations
 will be also equipped with a `classical' term $\displaystyle
 -\frac{d}{dt}\frac{\pd L}{\pd \dot{u}}$ on the left hand side. It
 can be shown that in this case infinitesimal criterion and
 conservation law can be obtained by combining the results of the
 present analysis (see Theorems \ref{th:inf criterion bez a},
 \ref{th:inf criterion} and \ref{th:fracNoether}) and the classical
 theory (see e.g. \cite{Olv}).
 \erem

 \section{Infinitesimal invariance}
 \label{sec:inf inv}

 Let $G$ be a local one-parameter group of transformations acting
 on a space of independent and dependent variables as
 follows: $(\bar{t},\bar{u})=g_\eta\cdot (t,u)=(\Xi_\eta (t,u),
 \Psi_\eta (t,u))$, for smooth functions $\Xi_\eta$ and
 $\Psi_\eta$, and $g_\eta\in G$.
 Let
 $$
 \vf=\tau(t,u)\frac{\pd}{\pd t}+\xi(t,u)\frac{\pd}{\pd u}
 $$
 be the infinitesimal generator of $G$. Then we also have
 \beq \label{eq:transformation group}
 \begin{array}{rcl}
 \bar{t}&=& t+\eta\tau(t,u)+o(\eta)\\[3pt]
 \bar{u}&=& u+\eta\xi(t,u)+o(\eta).
 \end{array}
 \eeq

 We introduce the following notation (cf., e.g., \cite{Vujanovic}):
 $
 \Delta t=\frac{d}{d\eta}|_{_{\eta=0}}(\bar{t}-t)
 $
 and
 $
 \Delta u=\frac{d}{d\eta}|_{_{\eta=0}}(\bar{u}(\bar{t})-u(t))
 $.
 More precisely, the notations $\Delta t$ and $\Delta u$ denote
 $\lim_{\eta\to 0}\frac{\bar{t}(\eta)-t}{\eta}$ and
 $\lim_{\eta\to 0}\frac{\bar{u}(\bar{t},\eta)-u(t)}{\eta}$,
 respectively.
 It follows from (\ref{eq:transformation group}) that
 $\Delta t=\tau$,
 and writing the Taylor expansion of $\bar{u}(\bar{t})=
 \bar{u}(t+\eta\tau(t,u)+o(\eta))$
 at $\eta=0$ yields that $\Delta u=\xi$.

 On the other hand, if we write Taylor expansion of
 $\bar{u}(\bar{t})$ at $\bar{t}=t$ we obtain
 $$
 \Delta u=\frac{d}{d\eta}{\Big\vert}_{_{\eta=0}}(\bar{u}(t)-u(t))+\dot{u}\Delta t,
 $$
 or, if we introduce the Lagrangian variation
 \beq \label{eq:Lagrangian variation}
 \delta u:=\frac{d}{d\eta}{\Big\vert}_{_{\eta=0}}(\bar{u}(t)-u(t))
 \eeq
 we obtain
 $$
 \Delta u=\delta u +\dot{u}\Delta t.
 $$
 Thus, we have
 $$
 \delta u=\xi - \tau\dot{u}.
 $$
 Note that $\delta t=0$.

 In the same way we can define $\Delta F$ and $\delta F$ of an
 arbitrary absolutely continuous function $F=F(t,u(t),\dot{u}(t))$:
 $$
 \Delta F=\frac{d}{d\eta}{\Big\vert}_{_{\eta=0}}\Big(
 F(\bar{t},\bar{u}(\bar{t}),\dot{\bar{u}}(\bar{t}))
 -F(t,u(t),\dot{u}(t)) \Big)=\frac{\pd F}{\pd t}\Delta t+
 \frac{\pd F}{\pd u}\Delta u +\frac{\pd F}{\pd \dot{u}}\Delta \dot{u}
 $$
 $$
 \delta F=\frac{d}{d\eta}{\Big\vert}_{_{\eta=0}} \Big(
 F(t,\bar{u}(t),\dot{\bar{u}}(t)) -F(t,u(t),\dot{u}(t))
 \Big)=\frac{\pd F}{\pd u}
 \delta u +\frac{\pd F}{\pd \dot{u}}\delta \dot{u},
 $$
 and
 $$
 \Delta F=\delta F + \dot{F}\Delta t.
 $$

 It will be apparent in the forthcoming sections that, if we
 want to find an infinitesimal criterion, we need to
 know $\Delta {}_aD_t^\alpha u$ and $\Delta \cL$.
 Therefore, we have to transform the
 left Riemann-Liouville fractional derivative of $u$ under
 the action of a local one-parameter group of transformations
 (\ref{eq:transformation group}).
 There are two different cases which will be considered
 separately. In the first one, the lower bound
 $a$ in ${}_aD_t^\alpha u$ is not
 transformed, while in the second case $a$ is transformed
 in the same way as the
 independent variable $t$. Physically, the
 first case is important when ${}_aD_t^\alpha u$ represents
 memory effects, and the second one is important when action on
 a distance is involved.

 \subsection{The case when $a$ in ${}_aD_t^\alpha u$ is not transformed}
 \label{subsec:a not transformed}

 In this section
 we consider a local group of transformations $G$ which transforms
 $t\in(A,B)$ into $\bar{t}\in(\bar{A},\bar{B})$ so that both
 intervals remain subintervals of $(a,b)$, but the action of $G$
 has no effect on the lower bound $a$ in ${}_aD_t^\alpha u$, i.e.,
 $\tau(a,u(a))=0$.
 So, suppose that $G$ acts on $t, u$ and ${}_aD_t^\alpha u$ in the
 following way:
 $$
 g_\eta\cdot (t, u, {}_aD_t^\alpha u):=(\bar{t}, \bar{u},
 {}_aD_{\bar{t}}^\alpha \bar{u}),
 $$
 where $\bar{t}$ and $\bar{u}$ are defined by
 (\ref{eq:transformation group}). In this case we have:

 \blem
 Let $u\in AC([a,b])$ and let $G$ be a local
 one-parameter group of transformations given by
 (\ref{eq:transformation group}). Then
 $$
 \Delta {}_aD_t^\alpha u =  {}_aD_t^\alpha \delta u +
 \frac{d}{dt} {}_aD_t^\alpha u\cdot \tau (t,u(t)),
 $$
 where
 $$
 \delta {}_aD_t^\alpha u = \frac{d}{d\eta}{\Big\vert}_{_{\eta=0}}\Big(
 {}_{a}D_{t}^\alpha \bar{u}- {}_aD_t^\alpha u \Big)
 = {}_aD_t^\alpha \delta u.
 $$
 \elem

 \pr
 To prove that $\delta {}_aD_t^\alpha u ={}_aD_t^\alpha \delta u$
 it is enough to apply the definition of the Lagrangian variation
 (\ref{eq:Lagrangian variation}). Also by definition we have
 $$
 \Delta {}_aD_t^\alpha u = \frac{d}{d\eta}{\Big\vert}_{_{\eta=0}}\Big(
 {}_{a}D_{\bar{t}}^\alpha \bar{u}- {}_aD_t^\alpha u \Big).
 $$
 Thus,
 \beas
 \Delta {}_aD_t^\alpha u &=& \frac{1}{\Gamma (1-\alpha)}
 \frac{d}{d\eta}{\Big\vert}_{_{\eta=0}} \Bigg[\
 \frac{d}{d\bar{t}} \int_{a}^{\bar{t}}
 \frac{\bar{u}(\theta)}{(\bar{t}-\theta)^\alpha}\, d\theta
 \pm \frac{d}{dt} \int_a^t \frac{\bar{u}(\theta)}{(t-\theta)^\alpha}\,
 d\theta
 -\frac{d}{dt} \int_a^t \frac{u(\theta)}{(t-\theta)^\alpha}\, d\theta
 \ \Bigg] \\
 &=& {}_a D_t^\alpha \delta u
 + \frac{d}{d\eta}{\Big\vert}_{_{\eta=0}}\Big(
 {}_{a}D_{\bar{t}}^\alpha \bar{u}- {}_aD_t^\alpha \bar{u} \Big) \\
 &=& {}_a D_t^\alpha \delta u
 + \frac{d}{d\bar{t}} \frac{d\bar{t}}{d\eta}{\Big\vert}_{_{\eta=0}}\Big(
 {}_{a}D_{\bar{t}}^\alpha \bar{u}- {}_aD_t^\alpha \bar{u} \Big)\\
 &=& {}_a D_t^\alpha \delta u + \frac{d}{dt} {}_aD_t^\alpha u
 \cdot \tau (t,u(t)).
 \eeas
 \ep

 \blem
 Let $\cL[u]$ be a functional of the form
 $
 \cL[u]=\int_A^B L(t, u(t), {}_aD_{t}^\alpha u)\, dt,
 $
 where $u$ is an absolutely continuous function in $[a,b]$,
 $(A,B)\subseteq (a,b)$
 and $L$ satisfies (\ref{eq:conditions on L}).
 Let $G$ be a local one-parameter group of transformations given by
 (\ref{eq:transformation group}).
 Then
 $$
 \Delta \cL=\delta \cL+(L\Delta t){\Big\vert}_{_{A}}^{^{B}},
 $$
 where
 $
 \delta \cL= \int_a^b \delta L\, dt.
 $
 \elem

 \pr
 Again it is clear that $\delta \cL= \int_a^b \delta L\,
 dt$. For $\Delta \cL$ we have
 \beas
 \Delta \cL &=&
 \frac{d}{d\eta}{\Big\vert}_{_{\eta=0}} \bigg(\
 \int_{\bar{A}}^{\bar{B}} L(\bar{t},
 \bar{u}(\bar{t}), {}_{a}D_{\bar{t}}^\alpha \bar{u})\, d\bar{t}
 - \int_A^B L(t, u(t), {}_aD_{t}^\alpha u)\, dt
 \ \bigg)\\
 &=&
 \frac{d}{d\eta}{\Big\vert}_{_{\eta=0}} \bigg(\
 \int_{A}^{B} L(\bar{t},
 \bar{u}(\bar{t}), {}_{a}D_{\bar{t}}^\alpha \bar{u})\,
 (1+\eta \dot{\tau} (t,u(t)))\, dt
 - \int_A^B L(t, u(t), {}_aD_{t}^\alpha u)\, dt
 \ \bigg) \\
 &=&
 \frac{d}{d\eta}{\Big\vert}_{_{\eta=0}} \bigg(\
 \int_{A}^{B} L(\bar{t},
 \bar{u}(\bar{t}), {}_{a}D_{\bar{t}}^\alpha \bar{u})\, dt
 -
 \int_{A}^{B}L(t,u(t), {}_{a}D_{t}^\alpha u)\, dt
 \ \bigg)\\
 && \quad
 + \int_{A}^{B} L(t,u(t), {}_{a}D_{t}^\alpha u)\,
 \dot{\tau} (t,u(t))\, dt.
 \eeas
 This further yields
 \beas
 \Delta \cL &=&
 \int_A^B \Delta L(t,u(t), {}_{a}D_{t}^\alpha u)\, dt
 +
 \int_{A}^{B} L(t,u(t), {}_{a}D_{t}^\alpha u)\,
 \dot{\tau} (t,u(t))\, dt \\
 &=&
 \int_A^B \delta L(t,u(t), {}_{a}D_{t}^\alpha u)\, dt
 +
 \int_{A}^{B} \frac{d}{dt}L(t,u(t), {}_{a}D_{t}^\alpha u)\,
 \tau (t,u(t))\, dt\\
 && \quad
 +
 \int_{A}^{B} L(t,u(t), {}_{a}D_{t}^\alpha u)\,
 \dot{\tau} (t,u(t))\, dt\\
 &=&
 \int_A^B \delta L\, dt
 +
 (L \tau)  {\Big\vert}^{^{B}}_{_{A}}
 \eeas
 and the claim is proved.
 \ep

 We define a variational symmetry group of the fractional
 variational problem (\ref{eq:fractional vp}), which we call
 a fractional variational symmetry group:

 \bd \label{def:variational symm group}
 A local one-parameter group of transformations $G$
 (\ref{eq:transformation group}) is a variational symmetry
 group of the fractional
 variational problem (\ref{eq:fractional vp}) if the following
 conditions holds: for every
 $[A',B']\subset (A,B)$,
 $u=u(t)\in AC([A',B'])$ and $g_\eta\in G$ such that
 $\bar{u}(\bar{t})=g_\eta\cdot u(\bar{t})$ is
 in $AC([\bar{A}',\bar{B}'])$, we have
 \beq \label{eq:variational symm group}
 \int_{\bar{A}'}^{\bar{B}'} L(\bar{t},\bar{u}(\bar{t}),
 {}_{a}D_{\bar{t}}^\alpha \bar{u})\, d\bar{t} =
 \int_{A'}^{B'} L(t,u(t),{}_aD_t^\alpha u)\, dt.
 \eeq
 \ed

 We are now able to prove the following infinitesimal criterion:

 \bthm \label{th:inf criterion bez a}
 Let $\cL[u]$ be a fractional variational problem
 (\ref{eq:fractional vp}) and let $G$ be a local one-parameter
 transformation group (\ref{eq:transformation group})
 with the infinitesimal generator $\vf
 =\tau(t,u)\pd_t+\xi(t,u)\pd_u$.
 Then $G$ is a variational symmetry group of $\cL$ if and only if
 \beq \label{eq:inf criterion}
 \tau\frac{\pd L}{\pd t}
  + \xi\frac{\pd L}{\pd u}
  + \Big({}_aD_t^\alpha (\xi-\dot{u}\tau)
    +\Big(\frac{d}{dt}{}_aD_t^\alpha u\Big)\tau \Big)\frac{\pd L}{\pd {}_aD_t^\alpha u}
  + L\dot{\tau} = 0.
 \eeq
 \ethm

 \pr
 Suppose that $G$ is a variational symmetry group of $\cL$. Then
 (\ref{eq:variational symm group}) holds for all subintervals
 $(A',B')$ of $(A,B)$ with closure $[A',B']\subset(A,B)$.
 We have:
 \beas
 \Delta\cL &=& \int_{A'}^{B'} \delta L\, dt + (L\Delta t) {\Big\vert}_{_{A'}}^{^{B'}} \\
 &=& \int_{A'}^{B'} \bigg( \frac{\pd L}{\pd u} \delta u
  + \frac{\pd L}{\pd {}_aD_t^\alpha u} \delta ({}_aD_t^\alpha u)\bigg)\, dt
  + (L\Delta t) {\Big\vert}_{_{A'}}^{^{B'}}\\
 &=& \int_{A'}^{B'} \bigg( \frac{\pd L}{\pd u} (\Delta u -\dot{u}\Delta t)
  + \frac{\pd L}{\pd {}_aD_t^\alpha u} {}_aD_t^\alpha (\Delta u -\dot{u}\Delta t)\bigg)\, dt
  + (L\Delta t) {\Big\vert}_{_{A'}}^{^{B'}}\\
 &=& \int_{A'}^{B'} \bigg( \frac{\pd L}{\pd u} \Delta u
  + \frac{\pd L}{\pd {}_aD_t^\alpha u} {}_aD_t^\alpha \Delta u\bigg)\, dt\\
 &&\quad - \int_{A'}^{B'} \bigg( \frac{\pd L}{\pd u} (\dot{u}\Delta t)
  + \frac{\pd L}{\pd {}_aD_t^\alpha u} {}_aD_t^\alpha (\dot{u}\Delta t)
  \pm \frac{\pd L}{\pd t} \Delta t \bigg)\, dt
  + (L\Delta t) {\Big\vert}_{_{A'}}^{^{B'}}\\
 &=& \ast
 \eeas
 Applying the Leibnitz rule for $\frac{d}{dt}(L\Delta t)$ we replace
 $
 \frac{\pd L}{\pd t} \Delta t +\frac{\pd L}{\pd u} \dot{u}\Delta t
 $
 by
 $
 \frac{d}{dt}(L\Delta t)
 - \frac{\pd L}{\pd {}_aD_t^\alpha u} (\frac{d}{dt}{}_aD_t^\alpha u)\Delta t
 -L\Delta t^{(1)}.
 $
 Then
 \beas
 \ast &=& \int_{A'}^{B'} \bigg(\frac{\pd L}{\pd t}\Delta t
  + \frac{\pd L}{\pd u} \Delta u
  + \frac{\pd L}{\pd {}_aD_t^\alpha u} {}_aD_t^\alpha \Delta u\\
 &&\quad + \frac{\pd L}{\pd {}_aD_t^\alpha u}
 \Big(\Big(\frac{d}{dt}{}_aD_t^\alpha u\Big)\Delta t
    - {}_aD_t^\alpha (\dot{u}\Delta t)\Big)
  +L\Delta t^{(1)}
  \bigg)\, dt.
 \eeas
 Since $\Delta\cL$ has to be zero in all $(A',B')$ with
 $[A',B']\subset (A,B)$, the above integrand has also to be
 equal zero, i.e.,
 $$
 \tau\frac{\pd L}{\pd t}
  + \xi\frac{\pd L}{\pd u}
  + \Big({}_aD_t^\alpha (\xi-\dot{u}\tau)
    +\Big(\frac{d}{dt}{}_aD_t^\alpha u\Big)\tau \Big)
    \frac{\pd L}{\pd {}_aD_t^\alpha u}
  +L\dot{\tau} = 0.
 $$
 where we have used that $\Delta u=\xi$ and $\Delta t=\tau$.
 Hence, necessity of the statement is proved. To prove that
 condition (\ref{eq:inf criterion}) is also sufficient,
 we first realize that if (\ref{eq:inf criterion}) holds on every
 $[A',B']\subset (A,B)$, then $\Delta\cL=0$ in every
 $[A',B']\subset (A,B)$. Thus, integrating $\Delta\cL$ from
 $0$ to $\eta$ we obtain (\ref{eq:variational symm group}), for
 $\eta$ near the identity. The proof is now complete.
 \ep

 In the following example we calculate the transformation
 of the fractional derivative ${}_aD_t^\alpha u$ under the group
 of translations. In Subsection \ref{subsec:a transformed} we
 shall perform such calculation also in the case when $a$ in
 ${}_aD_t^\alpha u$ is transformed.

 \bex \label{ex:translation, scaling}
 Let $G$ be a local one-parameter translation group: $(\bar{t},\bar{u})
 =(t+\eta, u+\eta)$, with the infinitesimal generator
 $\vf=\pd_t+\pd_u$. Then $\bar{u}(\bar{t})=u(\bar{t}-\eta)+\eta$
 and by a straightforward calculation it can be shown that
 $$
 {}_{a}D_{\bar{t}}^\alpha \bar{u} = {}_aD_t^\alpha (u+\eta) +
 \frac{1}{\Gamma (1-\alpha)} \frac{d}{dt}
 \int_{a-\eta}^{a}\frac{u(s)+\eta}{(t-s)^{\alpha}}\, ds.
 $$
 \eex

 \subsection{The case when $a$ in ${}_aD_t^\alpha u$ is transformed}
 \label{subsec:a transformed}

 Now let the action of a local one-parameter group of transformations
 (\ref{eq:transformation group}) be of the form
 $$
 g_\eta\cdot (t, u, {}_aD_t^\alpha u):=(\bar{t}, \bar{u},
 {}_{\bar{a}}D_{\bar{t}}^\alpha \bar{u}).
 $$
 Thus, the one-parameter group also acts on $a$ and transforms
 it to $\bar{a}$, where $\bar{a}=g_\eta\cdot t|_{t=a}$.
 This will also influence the calculation of $\Delta
 {}_aD_t^\alpha u$ and $\Delta \cL$:

 \blem
 Let $u\in AC([a,b])$ and let $G$ be a local
 one-parameter group of transformations
 (\ref{eq:transformation group}). Then
 \beq \label{eq:Delta of fractional der}
 \Delta {}_aD_t^\alpha u =  {}_aD_t^\alpha \delta u +
 \frac{d}{dt} {}_aD_t^\alpha u\cdot \tau (t,u(t)) +
 \frac{\alpha}{\Gamma (1-\alpha)} \frac{u(a)}{(t-a)^{\alpha +1}}
 \tau (a,u(a)),
 \eeq
 where
 $$
 \delta {}_aD_t^\alpha u = \frac{d}{d\eta}{\Big\vert}_{_{\eta=0}}\Big(
 {}_{a}D_{t}^\alpha \bar{u}- {}_aD_t^\alpha u \Big)
 = {}_aD_t^\alpha \delta u.
 $$
 \elem

 \pr
 Again using (\ref{eq:Lagrangian variation}) we check that
 $\delta {}_aD_t^\alpha u ={}_aD_t^\alpha \delta u$.
 To prove (\ref{eq:Delta of fractional der}) we start with
 the definition of $\Delta {}_aD_t^\alpha u$:
 $$
 \Delta {}_aD_t^\alpha u = \frac{d}{d\eta}{\Big\vert}_{_{\eta=0}}\Big(
 {}_{\bar{a}}D_{\bar{t}}^\alpha \bar{u}- {}_aD_t^\alpha u \Big).
 $$
 Thus,
 \beas
 \Delta {}_aD_t^\alpha u &=& \frac{1}{\Gamma (1-\alpha)}
 \frac{d}{d\eta}{\Big\vert}_{_{\eta=0}} \Bigg[\
 \frac{d}{d\bar{t}} \int_{\bar{a}}^{\bar{t}}
 \frac{\bar{u}(\theta)}{(\bar{t}-\theta)^\alpha}\, d\theta
 \pm \frac{d}{dt} \int_a^t \frac{\bar{u}(\theta)}{(t-\theta)^\alpha}\,
 d\theta
 -\frac{d}{dt} \int_a^t \frac{u(\theta)}{(t-\theta)^\alpha}\, d\theta
 \ \Bigg] \\
 &=& {}_a D_t^\alpha \delta u + \frac{1}{\Gamma (1-\alpha)}
 \frac{d}{d\eta}{\Big\vert}_{_{\eta=0}} \Bigg[\
 \frac{d}{d\bar{t}} \int_{\bar{a}}^{a}
 \frac{\bar{u}(\theta)}{(\bar{t}-\theta)^\alpha}\, d\theta
 + \frac{d}{d\bar{t}} \int_{a}^{\bar{t}}
 \frac{\bar{u}(\theta)}{(\bar{t}-\theta)^\alpha}\, d\theta\\
 && \quad
 - \frac{d}{dt} \int_a^t \frac{\bar{u}(\theta)}{(t-\theta)^\alpha}\, d\theta
 \ \Bigg] \\
 &=& {}_a D_t^\alpha \delta u + \frac{d}{dt} {}_aD_t^\alpha u
 \cdot \tau (t,u(t)) + \frac{1}{\Gamma (1-\alpha)}
 \frac{d}{d\eta}{\Big\vert}_{_{\eta=0}}
 \frac{d}{d\bar{t}} \int_{\bar{a}}^{a}
 \frac{\bar{u}(\theta)}{(\bar{t}-\theta)^{\alpha}}\, d\theta.
 \eeas
 Note that in the last term we can interchange
 the order of integration and differentiation with respect to
 $\bar{t}$. This gives
 $$
 \frac{-\alpha}{\Gamma (1-\alpha)} \frac{d}{d\eta}{\Big\vert}_{_{\eta=0}}
 \int_{\bar{a}}^{a}
 \frac{\bar{u}(\theta)}{(\bar{t}-\theta)^\alpha}\, d\theta.
 $$
 If we differentiate this integral with respect to $\eta$ at
 $\eta =0$ we eventually obtain (\ref{eq:Delta of fractional der}).
 \ep

 \blem
 Let $\cL[u]$ be a functional of the form
 $
 \cL[u]=\int_A^B L(t, u(t), {}_aD_{t}^\alpha u)\, dt,
 $
 where $u$ is an absolutely continuous function in $(a,b)$,
 $(A,B)\subseteq (a,b)$ and $L$ satisfies
 (\ref{eq:conditions on L}).
 Let $G$ be a local one-parameter group of transformations given by
 (\ref{eq:transformation group}).
 Then
 \beq \label{eq:Delta cL}
 \Delta \cL=\delta \cL+(L\Delta t){\Big\vert}_{_{A}}^{^{B}}
 +\frac{\alpha}{\Gamma (1-\alpha)}\frac{u(a)}{(t-a)^{\alpha +1}}
 \tau (a,u(a))\int_A^B \frac{\pd L}{\pd {}_{a}D_{t}^\alpha u}\, dt,
 \eeq
 where
 $
 \delta \cL= \int_A^B \delta L\, dt.
 $
 \elem

 \pr
 Clearly, $\delta \cL= \int_a^b \delta L\, dt$.
 On the other hand we have
 \beas
 \Delta \cL &=&
 \frac{d}{d\eta}{\Big\vert}_{_{\eta=0}} \bigg(\
 \int_{\bar{A}}^{\bar{B}} L(\bar{t},
 \bar{u}(\bar{t}), {}_{\bar{a}}D_{\bar{t}}^\alpha \bar{u})\, d\bar{t}
 - \int_A^B L(t, u(t), {}_aD_{t}^\alpha u)\, dt
 \ \bigg)\\
 &=&
 \frac{d}{d\eta}{\Big\vert}_{_{\eta=0}} \bigg(\
 \int_{A}^{B} L(\bar{t},
 \bar{u}(\bar{t}), {}_{\bar{a}}D_{\bar{t}}^\alpha \bar{u})\,
 (1+\eta \dot{\tau} (t,u(t)))\, dt
 - \int_A^B L(t, u(t), {}_aD_{t}^\alpha u)\, dt
 \ \bigg)\\
 &=&
 \int_{A}^{B} \Delta L\, dt
 + \int_{A}^{B} L(t,u(t), {}_{a}D_{t}^\alpha u)\,
 \dot{\tau} (t,u(t))\, dt\\
 &=& \int_{A}^{B} \Big(
 \frac{\pd L}{\pd t}\Delta t+ \frac{\pd L}{\pd u}\Delta u
 +\frac{\pd L}{\pd {}_{a}D_{t}^\alpha u}\Delta {}_{a}D_{t}^\alpha u
 \Big)\, dt
 + \int_{A}^{B} L(t,u(t), {}_{a}D_{t}^\alpha u)\,
 \dot{\tau} (t,u(t))\, dt.
 \eeas
 We now apply (\ref{eq:Delta of fractional der}) to obtain
 \beas
 \Delta \cL &=&
 \int_{A}^{B} \bigg(
 \frac{\pd L}{\pd t}\Delta t+ \frac{\pd L}{\pd u}(\delta u + \dot{u}\Delta t)
 +\frac{\pd L}{\pd {}_{a}D_{t}^\alpha u}
 \Big(
 {}_aD_t^\alpha \delta u +
 \frac{d}{dt} {}_aD_t^\alpha u\cdot \tau (t,u(t))\\
 && \quad
 + \frac{\alpha}{\Gamma (1-\alpha)} \frac{u(a)}{(t-a)^{\alpha +1}}
 \tau (a,u(a))
 \Big)
 \bigg)\, dt
 + \int_{A}^{B} L(t,u(t), {}_{a}D_{t}^\alpha u)\,
 \dot{\tau} (t,u(t))\, dt\\
 &=&
 \int_A^B \delta L\, dt + (L\Delta t){\Big\vert}_{_{A}}^{^{B}}
 + \frac{\alpha}{\Gamma (1-\alpha)} \frac{u(a)}{(t-a)^{\alpha +1}}
 \tau (a,u(a))\int_A^B \frac{\pd L}{\pd {}_{a}D_{t}^\alpha u}\,
 dt,\\
 &=&
 \delta \cL+(L\Delta t){\Big\vert}_{_{A}}^{^{B}}
 +\frac{\alpha}{\Gamma (1-\alpha)}\frac{u(a)}{(t-a)^{\alpha +1}}
 \tau (a,u(a))\int_A^B \frac{\pd L}{\pd {}_{a}D_{t}^\alpha u}\, dt.
 \eeas
 \ep

 \brem \label{rem:u(a)=0}
 It should be emphasized that the last summand on the right-hand side
 of (\ref{eq:Delta of fractional der}), as well as of (\ref{eq:Delta cL}),
 appears only as a consequence of the fact that the action of a transformation
 group affects also $a$, on which the left Riemann-Liouville fractional
 derivative depends. In the case when $u(a)=0$ or
 $\alpha=1$ that term is equal to zero (since $\lim_{\alpha\to 1^{-}}
 \Gamma(1-\alpha)=\infty$)
 and $\Delta \cL=\delta \cL+(L\Delta
 t){\Big\vert}_{A}^{B}$.
 \erem

 Now we define a variational symmetry group of the fractional
 variational problem (\ref{eq:fractional vp}) as follows:

 \bd \label{def:variational symm group 2}
 A local one-parameter group of transformations $G$
 (\ref{eq:transformation group}) is a variational symmetry
 group of the fractional
 variational problem (\ref{eq:fractional vp}) if the following
 conditions holds: for every
 $[A',B']\subset (A,B)$,
 $u=u(t)\in AC([A',B'])$ and $g_\eta\in G$ such that
 $\bar{u}(\bar{t})=g_\eta\cdot u(\bar{t})$ is
 in $AC([\bar{A}',\bar{B}'])$, we have
 $$
 \int_{\bar{A}'}^{\bar{B}'} L(\bar{t},\bar{u}(\bar{t}),
 {}_{\bar{a}}D_{\bar{t}}^\alpha \bar{u})\, d\bar{t} =
 \int_{A'}^{B'} L(t,u(t),{}_aD_t^\alpha u)\, dt.
 $$
 \ed

 Recall that we are solving the fractional variational problem
 (\ref{eq:fractional vp}) among all absolutely continuous functions
 in $[a,b]$, which
 additionally satisfy $u(a)=0$. Therefore, the last term in both
 (\ref{eq:Delta of fractional der}) and (\ref{eq:Delta cL})
 vanishes, and the infinitesimal criterion reads the same as in
 Theorem \ref{th:inf criterion bez a}:

 \bthm \label{th:inf criterion}
 Let $\cL[u]$ be a fractional variational problem defined by
 (\ref{eq:fractional vp}) and let $G$ be a local one-parameter transformation
 group defined by (\ref{eq:transformation group})
 with the infinitesimal generator $\vf=\tau(t,u)\pd_t+\xi(t,u)\pd_u$.
 Assume that $u(a)=0$, for all admissible functions $u$.
 Then $G$ is a variational symmetry group of $\cL$ if and only if
 \beq \label{eq:inf criterion 2}
 \tau\frac{\pd L}{\pd t}
  + \xi\frac{\pd L}{\pd u}
  + \Big({}_aD_t^\alpha (\xi-\dot{u}\tau)
    +\Big(\frac{d}{dt}{}_aD_t^\alpha u\Big)\tau \Big)\frac{\pd L}{\pd {}_aD_t^\alpha u}
  + L\dot{\tau} = 0.
 \eeq
 \ethm

 \pr
 See the proof of Theorem \ref{th:inf criterion bez a}.
 \ep

 \bex \label{ex:translation, scaling a}
 Let $\vf=\pd_t+\pd_u$ be the infinitesimal generator of the
 translation group $(\bar{t},\bar{u})=(t+\eta, u+\eta)$. Then
 $\bar{u}(\bar{t})=u(\bar{t}-\eta)+\eta$ and
 $$
 {}_{\bar{a}}D_{\bar{t}}^\alpha \bar{u}
 = {}_aD_t^\alpha (u+\eta).
 $$
 \eex

 \section{N\"{o}ther's theorem}
 \label{sec:Noether}

 In the formulation of N\"{o}ther's theorem we shall need a
 generalization of the fractional integration by parts.

 \blem \label{l:rRL to lRL}
 Let $f,g\in AC([a,b])$. Then, for all $t\in[a,b]$ the following
 formula holds:
 \beq \label{eq:rRL to lRL}
 \begin{array}{l}
 \displaystyle
 \int_a^t f(s)\cdot {}_sD_b^\alpha g\, ds =
 \int_a^t {}_aD_s^\alpha f \cdot g(s)\, ds \\
 \qquad \displaystyle
 + \int_a^t f(s)\cdot \frac{1}{\Gamma (1-\alpha)}
 \bigg[
 \frac{g(b)}{(b-s)^\alpha} - \frac{g(t)}{(t-s)^\alpha}
 - \int_t^b \frac{\dot{g}(\sigma)}{(\sigma -s)^\alpha}\, d\sigma
 \bigg]\, ds.
 \end{array}
 \eeq
 \elem

 \pr
 In order to derive (\ref{eq:rRL to lRL}) we shall use the
 representation of the
 right (resp.\ left) Riemann-Liouville fractional derivative via
 the right (resp.\ left) Caputo fractional derivative
 (\ref{eq:relation left rl c}), as well as
 (\ref{eq:frac int by parts}):
 \beas
 \int_a^t f(s)\cdot {}_sD_b^\alpha g\, ds &=&
 \int_a^t f(s)\cdot \frac{1}{\Gamma (1-\alpha)}
 \bigg[
 \frac{g(b)}{(b-s)^\alpha}
 - \int_s^b \frac{\dot{g}(\sigma)}{(\sigma -s)^\alpha}\, d\sigma
 \bigg]\, ds\\
 &=& \int_a^t f(s)\cdot \frac{1}{\Gamma (1-\alpha)}
 \bigg[
 \frac{g(b)}{(b-s)^\alpha}
 - \int_s^t \frac{\dot{g}(\sigma)}{(\sigma -s)^\alpha}\, d\sigma\\
 && \quad
 - \int_t^b \frac{\dot{g}(\sigma)}{(\sigma -s)^\alpha}\, d\sigma
 \pm \frac{g(t)}{(t-s)^\alpha}
 \bigg]\, ds\\
 &=& \int_a^t \Bigg[f(s)\cdot {}_sD_t^\alpha g
 +f(s)\cdot \frac{1}{\Gamma (1-\alpha)}
 \bigg[
 \frac{g(b)}{(b-s)^\alpha}\\
 && \quad
 - \frac{g(t)}{(t-s)^\alpha}
 - \int_t^b \frac{\dot{g}(\sigma)}{(\sigma -s)^\alpha}\, d\sigma
 \bigg] \Bigg]\, ds \\
 &=& \int_a^t \Bigg[{}_aD_s^\alpha f\cdot g(s)
 +f(s)\cdot \frac{1}{\Gamma (1-\alpha)}
 \bigg[
 \frac{g(b)}{(b-s)^\alpha}\\
 && \quad
 - \frac{g(t)}{(t-s)^\alpha}
 - \int_t^b \frac{\dot{g}(\sigma)}{(\sigma -s)^\alpha}\, d\sigma
 \bigg] \Bigg]\, ds.
 \eeas
 \ep

 \brem
 If we put $t=b$ in (\ref{eq:rRL to lRL}), we obtain
 (\ref{eq:frac int by parts}).
 \erem

 The main goal of symmetry group analysis in the calculus
 of variations are the first integrals of Euler-Lagrange
 equations of a variational problem, that is a N\"other type result.

 Analogously to the classical case, an expression
 $\frac{d}{dt}P(t,u(t), {}_aD_t^\alpha u)=0$ is called
 a fractional first integral
 (or a fractional conservation law)
 for a fractional differential
 equation $F(t,u(t), {}_aD_t^\alpha u)=0$, if it vanishes along all
 solutions $u(t)$ of $F$.

 In the statement which is to follow, we prove a version of
 the fractional N\"{o}ther theorem. As we will show, a fractional
 conserved quantity will also contain integral terms, which is
 unavoidable, due to the presence of fractional derivatives.
 If $\alpha=1$ then (\ref{eq:frac cl}) reduces to the well-known
 classical conservation laws for a first-order variational problem.

 \bthm \label{th:fracNoether}
 (N\"{o}ther's theorem)
 Let $G$ be a local fractional variational symmetry group
 defined by (\ref{eq:transformation group}) of the
 fractional variational problem (\ref{eq:fractional vp}), and let
 $\vf=\tau(t,u(t))\pd_t+\xi(t,u(t))\pd_u$ be the infinitesimal
 generator of $G$. Then
 \beq \label{eq:frac cl}
 L\tau+\int_a^t\bigg({}_aD_s^\alpha (\xi-\dot{u}\tau)\frac{\pd L}{\pd {}_aD_s^\alpha u}
 -(\xi-\dot{u}\tau){}_sD_B^\alpha\Big(\frac{\pd L}{\pd {}_aD_s^\alpha u}\Big)
 \bigg)\, ds = \mbox{ const.}
 \eeq
 is a fractional first integral (or fractional conservation law)
 for the Euler-Lagrange equation (\ref{eq:EL AB}).

 The fractional conservation law (\ref{eq:frac cl}) can
 equivalently be written in the form
 \bea
 && L\tau - \int_a^t (\xi-\dot{u}\tau)\cdot \frac{1}{\Gamma (1-\alpha)}
 \bigg[
 \frac{\pd_3L(B,u(B), {}_aD_B^\alpha u)}{(B-s)^\alpha} -
 \frac{\pd_3L(t,u(t), {}_aD_t^\alpha u)}{(t-s)^\alpha}\nonumber\\
 && \qquad\qquad\qquad\qquad -
 \int_t^b \frac{\frac{d}{d\sigma}\pd_3 L(\sigma, u(\sigma),
 {}_aD_\sigma^\alpha u)}{(\sigma-s)^{\alpha}}
 \, d\sigma \bigg]\, ds = \mbox{ const.} \label{eq:frac cl-2}
 \eea
 \ethm

 \brem
 This theorem is valid in the case when the lower bound $a$ in
 ${}_aD_t^\alpha u$ is not transformed, as well as in the case when
 $a$ is transformed.
 In the second case we in addition have to suppose that $u(a)=0$,
 for all admissible functions $u$, which provides the same form of
 the infinitesimal criterion in both cases.
 \erem

 \pr
 We want to insert the Euler-Lagrange equation (\ref{eq:EL AB RL})
 into the infinitesimal criterion (\ref{eq:inf criterion}) or
 (\ref{eq:inf criterion 2}). Hence, we will write the latter
 in a suitable form:
 \beas
 0&=& \tau\frac{\pd L}{\pd t}
  + \xi\frac{\pd L}{\pd u}
  + \Big({}_aD_t^\alpha (\xi-\dot{u}\tau)
    +\Big(\frac{d}{dt}{}_aD_t^\alpha u\Big)\tau \Big)\frac{\pd L}{\pd {}_aD_t^\alpha u}
  + L\dot{\tau}
  \pm \dot{u}\tau \frac{\pd L}{\pd u}\\
 &=& \tau \bigg(\frac{\pd L}{\pd t}+\frac{\pd L}{\pd u}\cdot \dot{u}
 +\frac{\pd L}{\pd {}_aD_t^\alpha u}\cdot\frac{d}{dt}{}_aD_t^\alpha u\bigg)
 + (\xi-\dot{u}\tau)\frac{\pd L}{\pd u}
 +{}_aD_t^\alpha (\xi-\dot{u}\tau)\frac{\pd L}{\pd {}_aD_t^\alpha u}
 +L\dot{\tau}\\
 &=& \tau \frac{d}{dt} L + L\frac{d}{dt}\tau+ (\xi-\dot{u}\tau)\frac{\pd L}{\pd u}
  + {}_aD_t^\alpha (\xi-\dot{u}\tau)\frac{\pd L}{\pd {}_aD_t^\alpha u}
 \pm (\xi-\dot{u}\tau){}_tD_B^\alpha \Big(\frac{\pd L}{\pd {}_aD_t^\alpha u}\Big)\\
 &=& \frac{d}{dt}\big(L\tau\big) + (\xi-\dot{u}\tau) \bigg(\frac{\pd L}{\pd u}+
  {}_tD_B^\alpha \Big(\frac{\pd L}{\pd {}_aD_t^\alpha u}\Big)\bigg)\\
 && \quad
  +\frac{d}{dt} \int_a^t\bigg({}_aD_s^\alpha (\xi-\dot{u}\tau)\frac{\pd L}{\pd {}_aD_s^\alpha u}
 -(\xi-\dot{u}\tau){}_sD_B^\alpha\Big(\frac{\pd L}{\pd {}_aD_s^\alpha u}\Big)
 \bigg)\, ds.
 \eeas
 The middle term in the last equation is the Euler-Lagrange equation
 (\ref{eq:EL AB RL}) multiplied by $\xi-\dot{u}\tau$. Therefore,
 $$
 \frac{d}{dt}\bigg(L\tau + \int_a^t\Big({}_aD_s^\alpha (\xi-\dot{u}\tau)
 \frac{\pd L}{\pd {}_aD_s^\alpha u}
 -(\xi-\dot{u}\tau){}_sD_B^\alpha(\frac{\pd L}{\pd {}_aD_s^\alpha u})
 \Big)\, ds \bigg)=0
 $$
 for all solutions of the Euler-Lagrange equation, and
 (\ref{eq:frac cl}) holds.
 If we now apply Lemma \ref{l:rRL to lRL} we can transform
 (\ref{eq:frac cl}) into (\ref{eq:frac cl-2}).
 \ep

 In the next example we consider the case in which the Lagrangian of
 a fractional variational problem does not depend on $t$ explicitly.

 \bex \label{ex:translation inv}
 Consider a variational problem whose Lagrangian does not depend
 on $t$ explicitly:
 \beq \label{eq:vp bez t}
 \cL[u]=\int_A^B L(u(t), {}_aD_t^\alpha u)\, dt.
 \eeq
 As in the classical case, the one-parameter group of translations
 of time is a variational symmetry group of (\ref{eq:vp bez t}).
 Indeed, $\bar{t}=t+\eta, \bar{u}=u$ (hence $\bar{u}(\bar{t})
 =u(\bar{t}-\eta)$) and ${}_{\bar{a}}D_{\bar{t}}^\alpha \bar{u}
 ={}_aD_t^\alpha u$ (see Example \ref{ex:translation, scaling a}).
 Therefore,
 $$
 \int_{\bar{A}}^{\bar{B}} L(\bar{u}(\bar{t}),
 {}_{\bar{a}}D_{\bar{t}}^\alpha \bar{u})\, d\bar{t}
 = \int_A^B L(u(t), {}_aD_t^\alpha u)\, dt.
 $$
 This fact can be also confirmed by the infinitesimal criterion:
 $$
 \frac{\pd L}{\pd {}_aD_t^\alpha u} \Big(-{}_aD_t^\alpha \Big(\frac{d}{dt} u\Big)
    +\frac{d}{dt}{}_aD_t^\alpha u\Big)
 = 0,
 $$
 where the last equality holds since $u(a)=0$. Thus,
 $\frac{d}{dt}{}_aD_t^\alpha u={}_aD_t^\alpha \frac{d}{dt} u$
 (see Introduction).

 Using N\"{o}ther's theorem we may write a fractional conservation
 law which comes from this translation group:
 $$
 L+\int_a^t\bigg(-{}_aD_t^\alpha \dot{u}\cdot\frac{\pd L}{\pd {}_aD_t^\alpha u}
 +\dot{u}\cdot{}_tD_B^\alpha\Big(\frac{\pd L}{\pd {}_aD_t^\alpha u}\Big)
 \bigg)\, ds = \mbox{ const.}
 $$
 or equivalently
 \bea
 && L + \int_a^t \dot{u}(s)\cdot \frac{1}{\Gamma (1-\alpha)}
 \bigg[
 \frac{\pd_3L(u(B), {}_aD_B^\alpha u)}{(B-s)^\alpha} -
 \frac{\pd_3L(u(t), {}_aD_t^\alpha u)}{(t-s)^\alpha} \nonumber\\
 && \qquad\qquad\qquad -
 \int_t^b \frac{\frac{d}{d\sigma}\pd_3 L(u(\sigma),
 {}_aD_\sigma^\alpha u)}{(\sigma-s)^{\alpha}}
 \, d\sigma \bigg]\, ds = \mbox{ const.} \label{eq:ex(i)-cl}
 \eea
 Note that (\ref{eq:ex(i)-cl}) tends to
 $$
 L-\dot{u}(t)\cdot\frac{\pd L(u(t),\dot{u}(t))}{\pd\dot{u}} = \mbox{ const.}
 $$
 when $\alpha\to 1^{-}$, which is the classical energy integral (see
 \cite{Vujanovic}).
 \eex

 \bex
 As a concrete example related to Example \ref{ex:translation inv},
 let us consider a fractional oscillator for which the
 Lagrangian reads:
 \beq \label{eq:ex(iii)}
 L=\frac{1}{2} ({}_aD_t^\alpha u)^2 - \omega^2\frac{u^2}{2},
 \eeq
 where $\omega$ is a constant (frequency).
 The fractional
 variational problem consists of minimizing the functional
 $\int_0^1 L\, dt$, among all smooth functions which satisfy
 the initial conditions $u(0)=0$ and $u^{(1)}(0)=1$.
 The Euler-Lagrange equation for such an $L$ is:
 \beq \label{eq:ex(iii)-EL}
 \omega^2 u- {}_tD_1^\alpha {}_0D_t^\alpha u=0.
 \eeq
 Since $L$ does not depend on $t$ explicitly and $u(0)=0$, the
 translation group $(\bar{t},\bar{u})=(t+\eta,u)$ is a fractional
 variational symmetry group for (\ref{eq:ex(iii)}). It generates
 the following fractional conservation law for
 (\ref{eq:ex(iii)-EL}):
 $$
 \frac{1}{2} ({}_aD_t^\alpha u)^2 - \omega^2\frac{u^2}{2} +
 \int_0^t \Big(- {}_0D_s^\alpha u^{(1)}\cdot {}_0D_s^\alpha u
 + u^{(1)}\cdot{}_sD_1^\alpha ({}_0D_s^\alpha u) \Big)\, ds = \mbox{const.}
 $$
 \eex

 \section{Approximations}
 \label{sec:approximations}

 In this section we use approximations
 of the Riemann-Liouville
 fractional derivatives by finite sums of derivatives of integer
 order, which leads to variational problems
 involving only classical derivatives. We examine the relation
 between the Euler-Lagrange equations, infinitesimal criterion
 and N\"{o}ther's theorem obtained in the process of
 approximation and the fractional Euler-Lagrange equations,
 infinitesimal criterion and N\"{o}ther's theorem
 derived in previous sections.

 In the sequel we make the following assumptions:
 \begin{itemize}
 \item[(i)] $L\in\cC^N([a,b]\times\R\times\R)$, at least, for some $N\in\N$.
 \item[(ii)] We will simplify the following calculations by considering the
 case $A=a$ and $B=b$.
 \item[(iii)] Let $(c,d)$, $-\infty<c<d<+\infty$, be an open interval in $\R$
 which contains $[a,b]$, such that for each $t\in [a,b]$ the closed
 ball $B(t, b-a)$, with center at $t$ and radius $b-a$, lies in $(c,d)$.
 \end{itemize}

 In addition we shall separately consider two cases:
 \begin{itemize}
 \item[(a)] Let $u\in\cC^\infty([a,b])$ such that $u(a)=a_0$,
 $u(b)=b_0$, for fixed $a_0,b_0\in\R$, and $L_3$ (where $L_3$
 stands for $\pd_3 L$) be a function
 in $[a,b]$ defined by $t\mapsto L_3(t)= L_3(t,u(t),{}_aD_t^\alpha u)$, $t\in
 [a,b]$. Let $L_3^{(i)}(b,b_0,p)=0$, for all $i\in\N$, meaning that for $(t,s,p)
 \mapsto L_3(t,s,p)$, $t\in [a,b]$, $s,p\in\R$, the following holds:
 $$
 \frac{\pd^i L_3}{\pd t^i}(b,b_0,p)=0;
  \quad \frac{\pd^i L_3}{\pd s^i}(b,b_0,p)=0;
  \quad \frac{\pd^i L_3}{\pd p^i}(b,b_0,p)=0, \quad \forall\,p\in\R.
 $$
 \item[(b)] Let $u\in\cC^\infty([a,b])$ such that $u^{(i)}(b)=0$,
 for all $i\in\N_0$, and $u(a)=a_0$, for fixed
 $a_0\in\R$. Let $L_3^{(i)}(b)=L_3^{(i)}(b,0,
 {}_aD_b^\alpha u)=0$, for all $i\in\N$ and for every fixed $u$,
 meaning that for $(t,s,p) \mapsto L_3(t,s,p)$, $t\in [a,b]$, $s,p\in\R$,
 the following holds:
 $$
 \frac{\pd^i L_3}{\pd t^i}(b,0,p)=0;
 \quad \frac{\pd^i L_3}{\pd p^i}(b,0,p)=0, \quad \forall\,p\in\R.
 $$
 \end{itemize}

 Let $f$ be a real analytic function in $(c,d)$. Then according to
 \cite[(15.4) and (1.48)]{SamkoKM} we have that
 \beq \label{eq:sum of derivatives}
 {}_aD_t^\alpha f =\sum_{i=0}^{\infty} {\alpha\choose i}
 \frac{(t-a)^{i-\alpha}}{\Gamma (i+1-\alpha)} f^{(i)}(t),
 \quad t\in B(t,b-a)\subset (c,d),
 \eeq
 where $\displaystyle {\alpha\choose i}=
 \frac{(-1)^{i-1}\alpha\Gamma(i-\alpha)}{\Gamma(1-\alpha)\Gamma(i+1)}$.

 Consider again the fractional variational problem
 (\ref{eq:fractional vp}). Assume that we are looking for a
 minimizer $u\in\cC^{2N}([a,b])$, for some $N\in\N$.
 We replace in the Lagrangian the left Riemann-Liouville fractional
 derivative ${}_aD_t^\alpha u$ by the finite sum of integer-valued
 derivatives as in (\ref{eq:sum of derivatives}):
 \beq \label{eq:N derivatives}
 \displaystyle
 \int_a^b L\Big(t, u(t), \sum_{i=0}^{N} {\alpha\choose i}
 \frac{(t-a)^{i-\alpha}}{\Gamma (i+1-\alpha)} u^{(i)}(t))\Big)\, dt
 =\int_a^b \bar{L} (t,u(t),u^{(1)}(t),u^{(2)}(t),
 \ldots,u^{(N)}(t))\, dt.
 \eeq
 Now the Lagrangian $\bar{L}$ depends on
 $t, u$ and all (classical) derivatives of $u$ up to order $N$.
 Moreover, $\pd_3 \bar{L},\ldots,\pd_{N+2} \bar{L}\in
 \cC^{N-1}([a,b]\times\R\times\R)$, since $\pd_i\bar{L}=\pd_3 L
 {\alpha \choose i} \frac{(t-a)^{i-\alpha}}{\Gamma(i+1-\alpha)}$,
 $i=3,...,N+2$.
 For the problem (\ref{eq:N derivatives}) we can
 calculate the Euler-Lagrange equations, infinitesimal criterion
 and conservation laws (via N\"{o}ther's theorem).

 Now we will consider the fractional variational problem
 (\ref{eq:fractional vp}) in the case (a) and in the case (b).

 The main result of this section is stated in the following
 theorem.

 \bthm \label{th:approx Noether}
 Let $\cL[u]$ be a fractional variational problem
 (\ref{eq:fractional vp}) which is being
 solved in the case (a) or (b). Denote
 by $CL$ the fractional conservation law (\ref{eq:frac cl}),
 and by $CL_N$ the fractional conservation law for the
 Euler-Lagrange equation (\ref{eq:EL-N}) which
 corresponds to the variational problem (\ref{eq:N derivatives}).
 Then
 $$
 CL_N\to CL \mbox{ in the weak sense, as } N\to \infty.
 $$
 \ethm

 Convergence in the weak sense here means that as a test
 function space we use the space of real analytic functions
 as follows.
 (Note that analytic functions are chosen in order to provide
 convergence of the sum as in (\ref{eq:sum of derivatives}).)

 Let $\mathcal{A}((c,d))$ be the space of real analytic
 functions in $(c,d)$ with the family of seminorms
 $$
 p_{[m,n]}(\vphi):= \sup_{t\in [m,n]}|\vphi(t)|, \quad\vphi\in
 \mathcal{A}((c,d)),
 $$
 where $[m,n]$ are subintervals of $(c,d)$. Every function
 $f\in\cC([a,b])$, which we extend to be zero in $(c,d)\backslash
 [a,b]$, defines an element of the dual $\mathcal{A}'((c,d))$ via
 $$
 \vphi\mapsto\langle f,\vphi \rangle =\int_a^b f(t)\vphi(t)dt,
 \quad \vphi\in\cA((c,d)).
 $$

 Before we prove the main Theorem \ref{th:approx Noether}, we shall
 prove several auxiliary results.

 First we recall a result from \cite{fractionalEL} which
 provides an expression for the
 right Riemann-Liouville fractional derivative in terms of the
 lower bound $a$, which figures in the left Riemann-Liouville
 fractional derivative.

 \bprop \label{prop:rrl-weak}
 Let $F\in\cC^\infty([a,b])$, such that $F^{(i)}(b)=0$,
 for all $i\in\N_0$, and $F\equiv 0$ in $(c,d)\backslash [a,b]$.
 Let ${}_tD_b^\alpha F$ be extended by zero in $(c,d)
 \backslash [a,b]$.
 Then:
 \begin{itemize}
 \item[(i)] For every $i\in\N$, the $(i-1)$-th derivative of
 $
 t\mapsto F(t)(t-a)^{i-\alpha}
 $
 is continuous at $t=a$ and $t=b$ and  the $i$-th derivative
 of this function, for $i\in\N_0$, is integrable in $(c,d)$ and
 supported in $[a,b]$.
 \item[(ii)] The partial sums $S_N,
 N\in\N_0$,
 $$
 t\mapsto S_N(t):=\left\{\begin{array}{ll}
 \displaystyle\sum_{i=0}^{N}
 \Big(-\frac{d}{dt}\Big)^i \bigg(F\cdot {\alpha\choose i}
 \frac{(t-a)^{i-\alpha}}{\Gamma (i+1-\alpha)}\bigg), & t\in [a,b],\\
 0, & t\in (c,d)\backslash [a,b],
 \end{array} \right.
 $$
 are  integrable functions in $(c,d)$ supported in $[a,b];$
 \item[(iii)]
 $$
 {}_tD_b^\alpha F = \sum_{i=0}^{\infty} \Big(-\frac{d}{dt}\Big)^i
 \bigg(F\cdot {\alpha\choose i}
 \frac{(t-a)^{i-\alpha}}{\Gamma (i+1-\alpha)}
 \bigg)
 $$
 in the weak sense.
 \end{itemize}
 \eprop

 \pr
 See \cite[Prop.\ 4.2]{fractionalEL}.
 \ep

 The Euler-Lagrange equation for (\ref{eq:N derivatives}) is of the
 following form:
 $$
 \sum_{i=0}^{N} \Big(-\frac{d}{dt}\Big)^i \frac{\pd \bar{L}}{\pd
 u^{(i)}}=0.
 $$
 This is equivalent to
 \beq \label{eq:EL-N}
 \frac{\pd L}{\pd u}+\sum_{i=0}^{N} \Big(-\frac{d}{dt}\Big)^i
 \bigg(\pd_3 L\cdot  {\alpha\choose i}
 \frac{(t-a)^{i-\alpha}}{\Gamma (i+1-\alpha)}\bigg)=0.
 \eeq

 \brem
 The Euler-Lagrange equation (\ref{eq:EL-N}) provides
 a necessary condition when one solves the variational problem
 (\ref{eq:N derivatives}) in the class $\cC^{2N}([a,b])$,
 with the prescribed boundary conditions at $a$ and $b$,
 i.e., $u(a)=a_0$ and $u(b)=b_0$, $a_0,b_0$ are fixed real numbers.
 \erem

 The following theorem shows that the Euler-Lagrange
 equation (\ref{eq:EL-N}) converges to (\ref{eq:EL popravljene}),
 as $N\to +\infty$, in the weak sense. To shorten the notation,
 we introduce $P_N$ and $P$ for the Euler-Lagrange equations in
 (\ref{eq:EL-N}) and (\ref{eq:EL popravljene}) respectively.

 \bthm \label{th:EL weak approx}
 Let $\cL[u]$ be the fractional variational problem
 (\ref{eq:fractional vp}) to be solved in the case (a) or (b).
 Denote by $P$ the fractional Euler-Lagrange equations
 (\ref{eq:EL AB}), and by $P_N$ the Euler-Lagrange equations
 (\ref{eq:EL-N}), which correspond to the variational problem
 (\ref{eq:N derivatives}), in which the left Riemann-Liouville
 fractional derivative is approximated according to
 (\ref{eq:sum of derivatives}) by a finite sum. Then in
 both cases a) and b)
 $$
 P_N\to P \mbox{ in the weak sense, as } N\to 0.
 $$
 \ethm

 \pr
 See \cite[Th.\ 4.3]{fractionalEL}.
 \ep

 Let
 \beq \label{eq:inf generator-N}
 \vf_N = \tau_N\frac{\pd}{\pd t}+\xi_N\frac{\pd}{\pd u}
 \eeq
 be the infinitesimal generator of a local one-parameter variational
 symmetry group of (\ref{eq:N derivatives}). The
 $N$-th prolongation of $\vf_N$ is given by
 $$
 {\rm pr}^{(N)} \vf_N = \vf_N+\sum_{i=1}^{N}\xi_N^i\frac{\pd}{\pd u^{(i)}},
 $$
 with
 $$
 \xi_N^i = \frac{d^i}{dt^i} \Big(\xi_N-\tau_N u^{(1)}\Big) +\tau_N
 u^{(i+1)}.
 $$

 \brem
 It should be emphasized that
 the vector field $\vf_N$ in (\ref{eq:inf generator-N}) differs (in
 general) from the one introduced at the beginning of Section
 \ref{sec:inf inv}.
 \erem

 \bthm \label{th:approx IC}
 Let $\cL[u]$ and $\bar{\cL}[u]$ be fractional variational
 problems (\ref{eq:fractional vp}) and (\ref{eq:N derivatives})
 respectively. Denote by $IC$ and $IC_N$ the corresponding
 infinitesimal criteria. If $\tau_N\to\tau$ and $\xi_N\to\xi$
 as $N\to\infty$, uniformly on compact sets, then
 $$
 IC_{N}\to IC \mbox{ uniformly on compact sets, as } N\to \infty.
 $$
 \ethm

 \pr
 The infinitesimal criterion for the variational problem
 (\ref{eq:N derivatives}) says that a vector field $\vf_N$
 generates a local one-parameter variational symmetry group of the
 functional (\ref{eq:N derivatives}) if and only if
 $$
 {\rm pr}^{(N)} \vf_N (\bar{L}) + \bar{L}\dot{\tau}_N =0.
 $$
 This is equivalent to
 \bea
 0 &=& \tau_N\frac{\pd\bar{L}}{\pd t}+\xi_N\frac{\pd\bar{L}}{\pd u}
 +\sum_{i=1}^{N}
 \bigg(\frac{d^i}{dt^i} \Big(\xi_N-\tau_N u^{(1)}\Big) +\tau_N
 u^{(i+1)}\bigg)
 \frac{\pd\bar{L}}{\pd u^{(i)}}
 + \bar{L}\dot{\tau}_N \nonumber\\
 &=& \tau_N\Big(\pd_1 L+\pd_3 L\sum_{i=0}^{N}{\alpha\choose i}
 \frac{(i-\alpha)(t-a)^{i-\alpha-1}}{\Gamma (i+1-\alpha)} u^{(i)}
 \Big)\nonumber\\
 && \quad
 +\xi_N\Big(\pd_2 L+\pd_3 L{\alpha\choose 0}
 \frac{(t-a)^{-\alpha}}{\Gamma (1-\alpha)}\Big) \nonumber\\
 && \quad
 +\sum_{i=1}^{N}
 \bigg(\frac{d^i}{dt^i} \Big(\xi_N-\tau_N u^{(1)}\Big) +\tau_N
 u^{(i+1)}\bigg)
 \pd_3 L {\alpha\choose i}\frac{(t-a)^{i-\alpha}}{\Gamma (i+1-\alpha)}
 + L\dot{\tau}_N \nonumber\\
 &=& \tau_N\pd_1 L+\xi_N\pd_2 Lv\nonumber\\
 && \quad
 +\Bigg[
 \sum_{i=1}^{N}{\alpha\choose i}\frac{(t-a)^{i-\alpha}}{\Gamma (i+1-\alpha)}
 (\xi_N-\tau_Nu^{(1)})^{(i)}
 + {\alpha\choose 0}\frac{(t-a)^{-\alpha}}{\Gamma (1-\alpha)}\xi_N
 \Bigg] \pd_3L \nonumber\\
 && \quad
 \pm {\alpha\choose 0}\frac{(t-a)^{-\alpha}}{\Gamma (1-\alpha)}\tau_Nu^{(1)}
 + \tau_N \Bigg[
 \sum_{i=0}^{N}{\alpha\choose i}
 \frac{(i-\alpha)(t-a)^{i-\alpha-1}}{\Gamma (i+1-\alpha)} u^{(i)} \nonumber\\
 && \quad
 +\sum_{i=1}^{N}{\alpha\choose i}
 \frac{(t-a)^{i-\alpha}}{\Gamma (i+1-\alpha)} u^{(i+1)}
 \Bigg]  \pd_3L
 + L\dot{\tau}_N \nonumber\\
 &=& \tau_N\pd_1 L+\xi_N\pd_2 L+
 \Bigg[
 \sum_{i=0}^{N}{\alpha\choose i}\frac{(t-a)^{i-\alpha}}{\Gamma (i+1-\alpha)}
 (\xi_N-\tau_Nu^{(1)})^{(i)} \nonumber\\
 && \quad
 + \frac{d}{dt}\bigg(
 \sum_{i=0}^{N}{\alpha\choose i}\frac{(t-a)^{i-\alpha}}{\Gamma (i+1-\alpha)}
 u^{(i)} \bigg) \tau_N
 \Bigg]\pd_3L
 + L\dot{\tau}_N \label{eq:inf crit-N-wo}
 \eea
 So, if $\tau_N\to\tau$ and $\xi_N\to\xi$ as $N\to\infty$, then
 the last expression on the right hand side tends to
 $$
 \tau\pd_1L + \xi\pd_2L
 + \Big({}_aD_t^\alpha (\xi-\dot{u}\tau)
 +\Big(\frac{d}{dt}{}_aD_t^\alpha u\Big)\tau \Big)\pd_3L
 + L\dot{\tau},
 $$
 and that is exactly the infinitesimal criterion
 (\ref{eq:inf criterion}) for (\ref{eq:fractional vp}).
 \ep

 \noindent{\bf Proof of Theorem \ref{th:approx Noether}.}
 First, we will derive the form of a conservation law for the
 Euler-Lagrange equation (\ref{eq:EL-N}) of the
 approximated variational problem (\ref{eq:N derivatives}), which
 comes from a variational symmetry group with infinitesimal
 generator (\ref{eq:inf generator-N}). For that purpose we need to
 insert the Euler-Lagrange equations (\ref{eq:EL-N}) into the
 infinitesimal criterion (\ref{eq:inf crit-N-wo}):
 \beas
 0 &=& \tau_N\pd_1 L+\xi_N\pd_2 L+
 \Bigg[
 \sum_{i=0}^{N}{\alpha\choose i}\frac{(t-a)^{i-\alpha}}{\Gamma (i+1-\alpha)}
 (\xi_N-\tau_Nu^{(1)})^{(i)} \\
 && \qquad
 + \frac{d}{dt}\bigg(
 \sum_{i=0}^{N}{\alpha\choose i}\frac{(t-a)^{i-\alpha}}{\Gamma (i+1-\alpha)}
 u^{(i)} \bigg) \tau_N
 \Bigg]\pd_3L
 + L\dot{\tau}_N
 \pm \tau_N u^{(1)} \pd_2L\\
 &=& \tau_N \Bigg( \pd_1 L +u^{(1)}\pd_2L + \pd_3L
 \frac{d}{dt}\bigg(
 \sum_{i=0}^{N}{\alpha\choose i}\frac{(t-a)^{i-\alpha}}{\Gamma (i+1-\alpha)}
 u^{(i)} \bigg)
 \Bigg)+ L\dot{\tau}_N\\
 && \qquad +\big(\xi_N-\tau_Nu^{(1)}\big)\pd_2L
 + \sum_{i=0}^{N}{\alpha\choose i}\frac{(t-a)^{i-\alpha}}{\Gamma (i+1-\alpha)}
 (\xi_N-\tau_Nu^{(1)})^{(i)} \pd_3L \\
 && \qquad \pm (\xi_N-\tau_Nu^{(1)})
 \sum_{i=0}^{N} \Big(-\frac{d}{dt}\Big)^i
 \bigg({\alpha\choose i}
 \frac{(t-a)^{i-\alpha}}{\Gamma (i+1-\alpha)} \pd_3 L\bigg)\\
 &=& \frac{d}{dt}\big(L\tau_N\big)
 + (\xi_N-\tau_Nu^{(1)}) \Bigg(\pd_2L + \sum_{i=0}^{N} \Big(-\frac{d}{dt}\Big)^i
 \bigg({\alpha\choose i}
 \frac{(t-a)^{i-\alpha}}{\Gamma (i+1-\alpha)} \pd_3
 L\bigg)\Bigg)\\
 && \qquad + \frac{d}{dt}\int_a^t
 \Bigg[
 \sum_{i=0}^{N}{\alpha\choose i}\frac{(s-a)^{i-\alpha}}{\Gamma (i+1-\alpha)}
 (\xi_N-\tau_Nu^{(1)})^{(i)} \pd_3L\\
 && \qquad \quad - (\xi_N-\tau_Nu^{(1)}) \sum_{i=0}^{N} \Big(-\frac{d}{ds}\Big)^i
 \bigg({\alpha\choose i}
 \frac{(s-a)^{i-\alpha}}{\Gamma (i+1-\alpha)} \pd_3 L\bigg)
 \Bigg]\, ds
 \eeas
 We recognize that the expression in brackets in the middle term in
 this sum is the Euler-Lagrange equation (\ref{eq:EL-N}), and hence
 vanishes. Therefore, we obtain that the following quantity is
 conserved:
 \beq \label{eq:cl-N}
 \begin{array}{l}
 \displaystyle
 L\tau_N
 + \int_a^t
 \Bigg[
 \sum_{i=0}^{N}{\alpha\choose i}\frac{(s-a)^{i-\alpha}}{\Gamma (i+1-\alpha)}
 (\xi_N-\tau_Nu^{(1)})^{(i)} \pd_3L \\
 \qquad \displaystyle
 - (\xi_N-\tau_Nu^{(1)}) \sum_{i=0}^{N}
 \Big(-\frac{d}{ds}\Big)^i
 \bigg({\alpha\choose i}
 \frac{(s-a)^{i-\alpha}}{\Gamma (i+1-\alpha)} \pd_3 L\bigg)
 \Bigg]\, ds
 =\mbox{const.}
 \end{array}
 \eeq
 This of course is one way to write the first integral of the
 Euler-Lagrange equations (\ref{eq:EL-N}), which corresponds to
 the local one-parameter variational symmetry group generated by
 (\ref{eq:inf generator-N}).

 Applying Proposition \ref{prop:rrl-weak} (ii) we obtain that
 (\ref{eq:cl-N}) converges to (\ref{eq:frac cl}) in the weak
 sense, provided that $(\tau_N,\xi_N)\to (\tau,\xi)$, as
 $N\to\infty$.
 \ep

 \brem
 Numerical analysis of fractional differential equations by the
 use of Theorems \ref{th:approx Noether}, \ref{th:EL weak approx}
 and \ref{th:approx IC} is not established and here we pose an
 open problem related to the use of our approximation procedure in
 applications.

 \erem

 \section{Concluding remarks}

 In this paper we studied variational principle containing
 fractional order derivatives. Especially, we obtained
 a condition which the symmetry group
 for this principle has to satisfy. This symmetry group
 is used for the formulation of a conservation law which generalizes
 the classical N\"{o}ther's theorem. In our work we used left
 Riemann-Liouville derivatives in Lagrangian density only. An
 extension including both left and right fractional derivatives
 could be easily performed. Using a method proposed in our
 previous work (see \cite{fractionalEL}) we approximated
 fractional Euler-Lagrangian equation with finite systems of
 integer order differential equations. In this way we obtained
 invariance conditions and corresponding conservation laws for
 such systems so that these conservation laws converge to the fractional
 conservation laws when the number of integer order equations tends
 to infinity.


 \end{document}